\documentclass[a4paper,12pt]{amsart}
\usepackage{amsmath,amsthm,latexsym,amssymb}
\makeatletter
\def\th@definition{%
  \normalfont
  }
\def\th@plain{%
  \slshape 
}
\def\th@remark{%
  \normalfont 
  \thm@preskip\topsep
  \divide\thm@preskip\tw@
  \thm@postskip\thm@preskip
}
\makeatother

\newcommand{\R}{\mathbb R}
\newcommand{\esup}{\mathop{\rm ess\,sup}}

\newcommand{\medint}{-\kern  -,375cm\int}

\let\disp=\displaystyle

\def\sign{\mathop{\rm sign}}

\def\esup{\mathop{\rm ess\,sup}}
\def\einf{\mathop{\rm ess\,inf}}

\def\dt{{d\over dt}}

\def\[{\bigl(}
\def\]{\bigr)}

\def\Box{\vrule height 6pt depth 1pt width 4pt}

\newcommand{\D}{\nabla}
\newcommand{\beq}{\begin{equation}}
\newcommand{\eeq}{\end{equation}}
\newcommand{\beqn}{\begin{eqnarray}}
\newcommand{\eeqn}{\end{eqnarray}}

\newcommand{\N}{\hbox{I\kern-.2em\hbox{N}}}

\theoremstyle{plain}
\newtheorem{myth}{Theorem}[section]
\newtheorem{myprop}{Proposition}[section]
\newtheorem{mylem}{Lemma}[section]
\newtheorem{mycor}{Corollary}[section]

\theoremstyle{definition}

\theoremstyle{remark}
\newtheorem{myrem}{Remark}[section]

\numberwithin{equation}{section}

\allowdisplaybreaks[1]

\begin{document}
\title[]{On isoperimetric inequalities\\
  with respect to infinite measures}

\author{ F. Brock,  A. Mercaldo, M.R. Posteraro}





\thanks{}
\date{}

\begin{abstract}
We study  isoperimetric problems with respect to infinite measures on $\R ^n$.
In  the case of the measure $\mu$ defined
by $d\mu = e^{c|x|^2 }\, dx$, $c\geq 0$, we prove  that, among all sets with given
$\mu-$measure, the ball centered
at the
origin  has the smallest (weighted) $\mu-$perimeter. Our results are then
applied
to obtain Polya-Szeg\"o-type inequalities, Sobolev embeddings theorems and
a comparison result for
elliptic boundary value problems.

\bigskip

\medskip

\end{abstract}

\maketitle

\section{Introduction}

Consider an elliptic  boundary value problem of the following type,
\begin{equation}
\left\{\begin{array}{ll}
-\hbox{div}(\varphi (x) |\nabla u|^{p-2}\nabla u)=f\varphi (x) & \hbox{in }
\Omega
\\
u=0 & \hbox{on } \partial\Omega,\label{bvp}
\end{array}
\right.
\end{equation}
where $\Omega$ is an
open subset of $\R^N$, possibly unbounded, and $f$ belongs to suitable weighted Lebesgue spaces.
We are interested in sharp explicit a-priori bounds for the weak solution to (\ref{bvp}).
Such type of problem can be examined by symmetrization methods.
However, the presence of the weight function in the operator in (\ref{bvp})
does not allow us to use the classical approach via
Schwarz symmetrization
given e.g.  in  \cite{[Ta1]}, and \cite{AT}. This leads us to introduce an appropriate
symmetrization based on a weighted isoperimetric inequality which is related to the
structure of the operator. A similar  approach which is based on the isoperimetric inequality for
Gauss measure has been carried out by the authors in \cite{[BBMP1]},
(see also \cite{[DiBF]},  \cite{[BCM]}).

In this paper we study isoperimetric inequalities for {\sl infinite } measures,
together with properties of corresponding
weighted symmetrizations.

To be more precise, let $\mu $ be a measure on $\R ^n $ defined by
$$
d\mu = \varphi (x)\, dx ,
$$
where $\varphi $ is a positive continuous function. For any smooth
set $\Omega $, we denote by
$$
P_{\mu } (\Omega ) = \int\limits _{\partial \Omega } \varphi (x)\,
    d{\mathcal H} _{n-1} (x)
$$
the weighted perimeter of $\Omega $ w.r.t. $\mu $ and for any fixed
number $m>0$, we denote by $I_{\mu } (m)$ the {\sl isoperimetric
function}, that is
\begin{equation}
\label{im} I_{\mu } (m) := \inf \{ P_{\mu } (\Omega ):\, \Omega \ \mbox{smooth}
, \mu (\Omega
)=m\}.
\end{equation}
We are interested in finding
 {\sl isoperimetric sets\/}, that is
smooth sets which realize the infimum in  (\ref{im}).

 Such a problem has
been treated in various settings. For example, if $\mu$ is the
Lebesque measure on $\R ^n $, then the isoperimetric sets are the
balls, i.e. if $\varphi (x)\equiv 1$, then $I_{\mu } (m)=P_{\mu }
(B)$, for any ball $B$ in $\R ^n $ with $\mu (B)=m$ (see for
instance \cite{[Ste]}, \cite{ [Ta3]}).

Moreover if $\mu$ is the Gauss measure, then  the isoperimetric sets are the
half-spaces of $\R ^n $, i.e. if $\varphi (x) = \mbox{exp } (-c|x|^2
)$ for some $c>0$ and $m\in (0,\mu (\R ^n ))$, then
 $I_{\mu } (m) =
P_{\mu } (H) $, where $H$ is any Euclidean half-space with $\mu
(H)=m$ (see for instance \cite{[Bo]}, \cite{[Eh]}, \cite{[Le]}).

Isoperimetric inequalities and their connections with rearrangements
have  received  considerable interest in the last decades (see
e.g. \cite{[Ka]}, \cite{[BuZa]}, \cite{[Ta2]}, \cite{[RS1]},
\cite{[RS2]} and the references cited therein). In the  paper
\cite{[BBMP2]}, the authors
recently analyzed symmetrizations w.r.t. {\sl finite\/} measures on $\R ^n $.

Here we investigate  {\sl infinite\/} measures on $\R ^n $ together with
Steiner and Schwarz symmetrizations.
One of our results is the following: If   $\mu $ is a measure defined by
\begin{equation} \label{mu} d\mu =  \,\mbox{exp } (c|x|^2 ) dx ,
\quad \mbox{for some
}\> c>0,
\end{equation}
then the only isoperimetric sets are Euclidean balls which are
centered at the origin, i.e. one has
\begin{equation} \label{im1}I_{\mu } (m) =
P_{\mu } (B_R ),
\end{equation}
where $R$ is chosen in such a way that $\mu (B_R)=m $.
Alternatively
one can express this isoperimetric inequality by using the notion of
weighted Schwarz symmetrization  $U^{ \displaystyle \star} $ of a
set $U$, which is the Euclidean ball centered at the origin such
that $\mu (U) =\mu (U^{\displaystyle \star}).$ With this notation,
(\ref{im1}) is equivalent to
\begin{equation}
\label{pm}
P_{\mu } (U) \ge P_{\mu } (U^{\displaystyle \star}),
\quad \mbox{for any smooth set }
U,
\end{equation}
with equality if and only if $U= U ^{\displaystyle \star} $
(see Theorem \ref{th6}).

We note that a proof of
inequality (\ref{pm}), but without the equality case, was given
by Borell already in 1986, in an unpublished preprint (see \cite{[Borellpreprint]}).
Theorem 4.1 was the subject of an earlier preprint of the authors,
and it has been presented at several conferences since 2005. R
ecently Benguria and Linde in  \cite{[BenguriaLinde]} used this result to obtain
eigenvalue bounds for the Dirichlet Schr\"odinger operator.
\\
We emphasize that after having finished our paper we learned that Theorem 4.1
has been also proved independently of us by Rosales, Canete, Bayle and Morgan (\cite{[Ros]}).
\\
All mentioned proofs are based on the observation that Steiner
symmetrization does not increase the perimeter.
Note, our proof differs from the one given in \cite{[Ros]} in that it does not
make use of the smoothness of a minimizing set. Therefore we include it here, for the convenience of the readers.
It is performed in several steps. \\
First we study the
one-dimensional case. More precisely, we consider a measure on $\R$ given
by
$$
d\mu_1 = \psi (x_1)\, dx_1,
$$
where $\psi$ is an even, positive and continuous function on $\R$ and
we prove that (\ref{im1}) holds iff $\psi$  is log-convex.
Then we consider a more general measure $\mu$ whose density is the product of two continuous
functions $\psi: \R \rightarrow \R_+$, and $\rho: \R^{n-1}
\rightarrow \R_+$ depending on  $x_1$ and $x^\prime=(x_2,.....x_n)$ respectively, i.e.
$$
d\mu = \psi (x_1)\,\rho (x^\prime) \, dx;
$$
we prove that  Steiner  symmetrization  with respect to this measure
decreases the perimeter (see Section 3.2 for the definition of Steiner symmetrization).
The last step
consists in  approximating the ball $U^{\displaystyle \star}$ by an
appropriate sequence of consecutive Steiner symmetrizations. In view of the product structure of the
density and the invariance  w.r.t. rotations of the measure (\ref{mu}), this leads to (\ref {pm}).

Isoperimetric inequality (\ref{pm})
has various consequences. For example, by  Talenti's result (\cite{[Ta2]}),
inequality (\ref{pm}) implies a Polya-Szeg\"o type inequality (see Theorem  \ref{Polya})
and the equality case is also studied. Moreover, we prove a Sobolev type
imbedding theorem in a weighted space w.r.t. the measure $\mu $ defined in (\ref{mu})
(see Theorem \ref{sobolev} ). The best constant in such an inequality is
obtained in a special case (see Corollary \ref{cor3}).

We now outline the content of the paper. In Section  2 we introduce some notation
and provide basic information about the weighted perimeter.
In Section 3 we study the Steiner symmetrization w.r.t. general measures.
In particular we prove that   Steiner symmetrization decreases the
weighted perimeter of a set having given measure. We also show a general form of
Polya-Szeg\"o type inequalities principle (see Theorem  \ref{ps}).
The proof of inequality (\ref{pm})   is given in Section 4.
Finally, Section 5 contains the comparison result

\section{Notation and preliminaries}
\setcounter{section}{2}
\setcounter{equation}{0}
In the whole paper $\mu $ will denote a
measure on $\R ^n $ defined by
\begin{equation}
\label{2.1}
d\mu = \varphi (x)\, dx ,
\end{equation}
where $\varphi \in C(\R ^n )$ and $\varphi (x)>0$ for any
$ x\in \R ^n $. Moreover $\mu_1 $ will denote a
measure on $\R $ defined by
\begin{equation}
\label{2.1bis}
d\mu_1 = \psi (x)\, dx ,
\end{equation}
where $\psi \in C(\R   )$ and $\psi (x)>0$ for any
$ x\in \R  $.
\\
We will always assume that the measures are infinite, that is
\begin{equation}
\label{2.100}
\mu (\R ^n )=+\infty ,\quad \mu _1 (\R) =+\infty .
\end{equation}
By ${\mathcal G M } ^n $ we denote the set of $\mu$-measurable sets with
finite
Lebesgue measure and by
${\mathcal F} ^n $ the set of all $\mu$-measurable
functions on $\R ^n $ such that  $\{ x:u(x)>t\} \in {\mathcal M} ^n$ for
every
$t>\inf u$.
\\
Let $\Omega $ be a domain in $\R ^n $ and  $p\in [1,+\infty )$. We
denote by
$L^p (\varphi ,\Omega )$ the space of $\mu$-measurable functions $u$
such
that
\begin{equation}\label{lorn}
\Vert u\Vert _{p ,\Omega } \equiv \left( \int_{\Omega }
|u|^p \, d\mu \right)
^{1/p} <+\infty ,
\end{equation}
endowed with the norm  (\ref{lorn}).
Furthermore, let $W^{1,p} (\varphi ,\Omega ) $
denote the weighted Sobolev space
containing all functions $u\in L^p (\varphi ,\Omega ) $
with weak derivatives
$u_{x_i } \in L^p (\varphi ,\Omega )$, $i=1,\ldots ,n$, and let
\begin{equation}
\label{2.6}
|\Vert u\Vert | _{p,\Omega } := \Vert u\Vert _{p,\Omega } +
\Vert \nabla u\Vert _{p,\Omega }
\end{equation}
be a norm in this space.
Finally, let $W_0 ^{1,p} (\varphi ,\Omega )$ be
the closure of $C_0 ^{\infty } (\Omega )$ under the norm (\ref{2.6}).
If $\Omega =\R ^n $ in one of the above spaces,
then we will omit the subindex $\Omega $ in the norms.
\\
For subsets $A,B,M$ of $ \R ^n $, let
$ A+B = \{ x+y :\, x\in A ,y\in B\} $
denote the Minkowski sum of $A$ and $B$,
and
$$
M_r :=  \{ x\in \R ^n :\, \mbox{dist } \{ x; M\} <r\}
= M+B_r ,\quad r>0 ,
$$
the exterior parallel sets of $M $, where $B_r$ is the
ball of radius $r$ with the center at the origin. \\
We will call a set $M\subset \R ^n $ {\sl smooth\/}, if $M$ is bounded,
 open, and if there is a number
$\varepsilon >0 $ such that for every $x_0 \in \partial M$, $\partial
M\cap
B_{\varepsilon } (x_0 )$ is a Lipschitz graph,
and
$M\cap B_{\varepsilon } (x_0 )$ lies on one side of $\partial M\cap
B_{\varepsilon } (x_0 )$ only. Observe that, this definition
includes polyhedra and excludes the presence of ``veils" and inner
``slices"
of $M$.
\\
If $M$ is a Borel set then we denote by
$\mu ^+ (M)$ the (lower outer)
{\sl Minkowski $\mu $-content\/} of the boundary of $M$
(see, for instance, \cite{[BuZa]}, p.69)
which is defined by
\begin{equation}
\label{2.2}
\mu ^+ (M) = \liminf _{r\to 0^+ } \frac{ \mu (M_r) -\mu (M)}{r}.
\end{equation}
If $\varphi \in W^{1,1}_{loc} (\R ^n )$, then we define the {\sl
$\mu$-perimeter\/} (in the sense of De Giorgi) by
\begin{equation}
\label{2.3}
P_{\mu} (M) = \sup \Big\{ \int_M \mbox{div} \, (\varphi (x) h(x))\,
dx :\
h\in C_0 ^1 (\R ^n ,\R ^n ),\, |h|\leq 1 \} .
\end{equation}
The following properties are well-known for the Lebesgue measure (see,
for
instance, \cite{[BuZa]}) and their proofs carry over to  general
measures $\mu $ defined in (\ref{2.1}) .\\
{\bf 1)} Both
$\mu ^+ (M)$ and
$P_{\mu } (M)$ can be seen as a ``weighted surface measure" of $M$,
that is, if $M$ is a smooth set then
\begin{equation}
\label{2.4}
\mu ^+ (M) =P_{\mu} (M) = \int_{\partial M } \varphi (x)\, {\mathcal H}
_{n-1} (dx) <+\infty .
\end{equation}
{\bf 2)} {\sl Lower semicontinuity\/}: If $\{ M_k \} \subset {\mathcal M}
^n$,
$M\in {\mathcal M} ^n$, and $\mu (M\Delta M_k \} \to 0 $ as $k\to
\infty $, then
\begin{equation}
\label{prop1}
\liminf_{k\to \infty} P_{\mu } (M_k ) \geq P_{\mu } (M).
\end{equation}
{\bf 3)} If $M\in {\mathcal M} ^n$, and $P_{\mu }(M)<\infty $,
then there is a sequence of smooth sets
$\{ M_k \} $ such that $\mu (M\Delta M_k )\to 0$ as $k\to \infty $,
and such
that
\begin{equation}
\label{prop2}
\lim_{k\to \infty } P_{\mu } (M_k ) = P_{\mu } (M).
\end{equation}
{\bf 4)} If $M$ is a Borel set in $\R ^n $, then
\begin{equation}
\label{prop3}
P_{\mu } (M)\leq \mu ^+ (M).
\end{equation}
\bigskip

We mention that
the theory of sets with finite $\mu $-perimeter
is imbedded  in the framework  of  BV-functions space, $BV(\varphi ,\R ^n )$, defined
as the set of all functions
$u\in L^1 (\varphi , \R ^n )$ such that
\begin{equation}\label{2.400}
\Vert Du \Vert _{BV}   :=
\sup \Big\{ \int\limits_{\R^n   }
    u(x) \mbox{ div }\Big( \varphi (x) h (x)\Big)  \, dx ,\
h \in C_0 ^1 (\R^n  ,\R ^n ),\
|h| \leq 1\Big\}  <+\infty .
\end{equation}
Notice that if $M$ has finite $\mu $-perimeter then
    the characteristic function of $M$,
    $\chi _M$ belongs to
$BV(\varphi ,\R ^n )$ and moreover $\Vert D\chi _M \Vert _{BV} = P_
{\mu } (M)$.
Furthermore, if $u\in W^{1,1} (\varphi ,\R ^n  ) $
then $\Vert Du \Vert _{BV} = \Vert \nabla u \Vert _1 $.
  \bigskip

Finally  we recall some well-known derivation formulas (cf.
\cite{[ADLT]},  \cite{[FM]} and \cite{[CiFu]}).
\\
We set  $\R _0 ^+ = [0,+\infty )$ and
we denote a point $x\in \R ^n$ by
$x=(x_1 ,x' )$ where $x_1\in \R$ and $x'\in \R^{n-1}$ $(n\ge2)$.
\\
Let $u\in {\mathcal F } ^n$, and let $m_u $ denote its $\mu  $-distribution
function w.r.t. the variable $x_1 $, that is$^{1}$
\setcounter{footnote}{1}
\footnotetext{\scriptsize{Here and in the following we will write
    $\{ u(\cdot ,x' )>t \} $ for $\{ x_1:\, u(x_1 ,x' )>t\}$. }}
\begin{equation}
\label{5.2.5}
m_u (t,x' ) := \mu _1\left (\{ u(\cdot ,x' )>t \}\right ) ,
\quad \ t\in \R,~  x' \in \R ^{n-1}.
\end{equation}
We also set
$$
u_- :=  \mbox{ess inf }u,
$$
and
\begin{equation}\label{vu}
V_u  :=   \{ (t,x'):\, t>u_-, \, x' \in \R ^{n-1} \} .
\end{equation}
Let
$u\in W^{1,p} (\varphi ,\R ^n ) $ for some $p\in [1,+\infty )$. We
define
\begin{equation}\label{du}
D_u (x'):= \{ x_1\in \R :\, u_{x_1 } (x_1 ,x' )=0 \} ,\quad x' \in \R
^{n-1} ,
\end{equation}
and the superlevel sets
\begin{equation}\label{eu}
E_u (t,x' ):=
\{ u(x_1 ,x' )>t\} \equiv \{ u(\cdot ,x' )>t\} ,
\quad (t,x' )\in V_u .
\end{equation}
Notice that, since $u(\cdot ,x') \in W^{1,p} (\varphi ,\R )$ for
a.e. $ x' \in \R ^{n-1}$, the Sobolev imbedding theorem tells us
that $u (\cdot ,x' )$
    is continuous for a.e. $x' \in \R ^{n-1} $. Hence
$E_u (t,x')$  is open and
$\partial E_u (t,x' )$ is countable for a.e. $(t,x')\in V_u $.
\\
Let us consider $g\in L^1 (\R ^n )$. Defining
$$
F(t,x' ) := \int\limits_{E_u (t,x') } g(x_1 ,x' )\, dx_1 ,
\quad (t,x')\in \R \times \R ^{n-1} ,
$$
we then have that $F(\cdot ,x' ) \in BV (\R ) $ for a.e.
$x' \in \R ^{n-1} $, and the Fleming-Rishel formula
(see, for instance, \cite{[Fe]}) tells us that
\begin{equation}
\label{5.4.3}
\frac{\partial}{\partial t} F(t,x') = \int\limits_{\partial E_u (t,x')}
\frac{g(\cdot ,x' )}{|u_{x_1 } (\cdot ,x')|} \, d{\mathcal H} _0
\qquad \mbox{for a.e. $(t,x') \in V_u $.}
\end{equation}
Moreover, we have by the co-area formula (see, e.g., \cite{[Fe]})
\begin{equation}
\int\limits_{\R } F(t,x') \, dt =
\int\limits_{\R \setminus D_u (x') }  g(x_1 ,x' )\, dx_1 =
\int\limits_{u_- (x')} ^{u_+ (x')} \int\limits_{\partial E_u (t,x')}
\frac{g(\cdot ,x')}{|u_{x_1}
(\cdot ,x' )|} \, d{\mathcal H} _0 \, dt,
\end{equation}
for a.e. $x'\in \R ^{n-1} $.
\\
Now let us assume that $u$ is a smooth function and satisfies
\begin{equation}
\label{5.4.6}
{\mathcal L} ^1 (D_u (x' )) =0 \quad \mbox{ for a.e. $ x' \in
\R^{n-1} $.}
\end{equation}
Here ${\mathcal L}^1 $ denotes the $1$-dimensional Lebesgue measure.
Then the following derivation formula holds true  (cf \cite{[ADLT]};
see also \cite{[FM]} and \cite{[CiFu]}).
\begin{mylem}
If $g\in C^1$, then we have
\begin{equation} \label{5.4.7}
\frac{\partial }{ \partial x_i }
\int\limits_{E_u(t,x')} g(x_1 ,x')\, dx_1=
\int\limits_{
E_u (t,x')} \frac{\partial g}{ \partial x_i }  (x_1 ,x') \, dx_1 +
\int\limits_{
\partial E_u (t,x')} g(\cdot ,x')
    \frac{u_{x_i } (\cdot ,x') }{|u_{x_1 } (\cdot ,x')|}  \, d{\mathcal
H} _0.
\end{equation}
\end{mylem}
\bigskip

\section{Measure on the real line}

Let $\psi $ a positive, even and  continuous function on $\R$, and
such that
$$
\int\limits_{\R} \psi (x)\, dx =+\infty .
$$
We then define a measure $\mu _1 $ on $\R $ by
$$
d\mu _1 = \psi (x)dx.
$$
Hence the primitive of $\psi $,
$$
\Psi (x):= \int\limits_{0 } ^x \psi (t)\, dt ,\quad x \in \R ,
$$
is strictly increasing and odd.
We introduce a continuous, positive function $J $ by
\begin{equation}
\label{1.200}
J(y)= \psi (\Psi ^{-1} (y)), \quad y\in \R ,
\end{equation}
where $\Psi ^{-1} : \R \to \R $
is the inverse of $\Psi $. Notice that $J$ is even, and (\ref{1.200})
implies that
\begin{equation}
\label{1.201}
\Psi ^{-1} (y) = \int\limits_{0} ^y \frac{dt}{J (t)} ,
\quad y\in \R .
\end{equation}
If $M \in {\mathcal M} ^1 $, and if $\mu _1 (M)<+\infty $, then
there exists a unique number
$c\geq 0 $ such that
\begin{equation}
\label{1.4}
\mu _1 ((-c,c) ) = \mu _1 ( M),
\end{equation}
and we set
\begin{equation}
\label{1.5}
M^* = (-c , c ).
\end{equation}
If $M\in {\mathcal M} ^1 $, and if $\mu _1 (M)=+\infty $ then we set
$M^* =\R $.
We call the set $M^* $
the $\mu _1 $-{\sl symmetrization\/} of $M$.
\bigskip

\begin{myrem}
 By the above definition, $M^* $ is some centered interval
$(-c ,c)$, ($c\in \R _0 ^+ \cup \{ +\infty \} $).
Notice that since we deal with measurable sets we usually do not
distinguish between two sets $M,N$ which are {\sl equivalent\/}, that
is,
which satisfy
$\mu _1 (M\Delta N)=0$.
On the other hand, the above definition can be easily precised if
$M$ is open or closed:
\\
If $M$ is open then the above settings have to be understood in
{\sl pointwise\/}
sense, so that $M^* $ is then open, too. If $M$ is closed then we
replace the open intervals in the above definition by closed intervals,
so that
$M^* $ is closed.
\end{myrem}

\begin{myrem}
Obviously, by definition
\begin{equation}
\label{1.1000}
\mu _1( M) = \mu _1 (M^* ) \quad \forall M\in {\mathcal M} ^1 .
\end{equation}
It is also easy to confirm the following monotonicity properties,
($M,N\in {\mathcal M} ^1 $):
\begin{eqnarray}
\label{1.101}
     & & M\subset N \Longrightarrow M^* \subset N^* ,
\\
\label{1.102}
     & & M^* \cup N^* \subset (M\cup N )^* ,\quad
M^* \cap N^* \supset (M\cap N ) ^* ,
\\
\label{1.103}
     & & \mu _1 (M\setminus N) \geq \mu _1 (M^* \setminus N^* ),\quad
\mu _1 (M\Delta N) \geq \mu _1 (M^* \Delta N^* ).
\end{eqnarray}
\end{myrem}
\bigskip

We now ask for a condition on the measure $\mu _1 $ such that the
$\mu _1 $-rearrangement decreases the
perimeter, that is we ask for a condition such that
the following isoperimetric inequality holds
\begin{equation}
\label{1.9}
P_{\mu _1 }(M) \geq P_{\mu _1 } (M^* ) \equiv
J(\mu _1 (M)),\quad
\forall M\in {\mathcal M} ^1 .
\end{equation}
Such a condition is given by Theorem \ref{th1} below
\begin{myth}\label{th1}
Inequality (\ref{1.9}) holds iff $J $ is convex.
Furthermore, let us assume that equality holds in (\ref{1.9}); if $J$ is
convex,
then $M$ is equivalent to an
interval, while if
$J$ is strictly convex, then $M=M^* $.
\end{myth}
\bigskip

\begin{myrem}\label{rem3}
$J$ is convex iff $\log \psi $ is convex (or equivalently,
if $\psi$
is log-convex).
A typical case is
$$
\psi (x) = e^{c|x| ^2  } ,
$$
where $c\in \R _0 ^+ $. Obviously if $c>0$
then
$J$ is {\sl strictly\/} convex.
\end{myrem}
\medskip

\begin{myrem} \label{rem4}
Isoperimetric inequality  (\ref{1.9}) and  property
(\ref{prop3}) imply an isoperimetric inequality
for the Minkowsky $\mu_1-$content, i.e.
\begin{equation}
\label{inequMink}
\mu_1 ^+(M)\geq \mu _1 ^+ (M^*) \ \
\mbox{ for every Borel set $M$.}
\end{equation}
Moreover it is easy to see that this also implies
$$
\mu _1 (M_r )\geq \mu _1 ((M^* ) _r )\quad \forall r>0 ,
$$
which is equivalent to
\begin{equation}
\label{1.13}
(M_r ) ^* \supset (M^* )_r  \qquad \forall
r>0.
\end{equation}
\end{myrem}
\medskip

\noindent{\bf Proof of Theorem \ref{th1} :\/}
First let us assume that $\mu _1 $ satisfies (\ref{1.9}), and let $I$ be any
finite interval
$(a,b)$, ($a<b$). Setting $\alpha =\Psi (a)$, $\beta = \Psi (b) $, we
have
that $I^* =(-c,c)$,
where $\Psi (c) = (\beta -\alpha ) /2$, then  (\ref{1.9}) reads as
$$
J(\alpha ) +J(\beta ) \geq 2J\Big( \frac{\beta -\alpha }{2} \Big) .
$$
Since $J(\alpha )=J(-\alpha )$, this implies that
$$
J(s) +J(t) \geq 2J\Big( \frac{s+t}{2} \Big) \quad \forall s,t \in \R ,
$$
which means that
$J$ is convex.
\\
Now let us assume that  $J$ is convex, and let $M$ be a smooth set.
Then $\displaystyle M=
\displaystyle\cup _{i=1} ^m (a_i ,b_i )$ where $m\in \N $, $a_i <b_i $,
and the intervals $[a_i ,b_i ]$ are mutually disjoint, $i=1,2,
\ldots ,m$.
Setting
$\alpha _i = \Psi (a_i )$, $\beta _i = \Psi (b_i )$, $i=1,\ldots ,m$,
we find, using the
convexity and evenness of $J$,
\begin{eqnarray*}
\qquad \ \ \ \ \ \ \ P_{\mu _1 } (M) & = &
\sum\limits_{i=1} ^m (J(\alpha _i )+J(\beta _i ))  \geq
\sum\limits_{i=1} ^m 2J\Big( \frac{\beta _i-\alpha _i }{2} \Big)
\\
    & \geq &
2m J\Big( \sum\limits_{i=1} ^m \frac{\beta _i -\alpha _i }{2m} \Big)
\geq
2 J
\Big( \sum\limits_{i=1} ^m \frac{\beta _i -\alpha _i }{2} \Big) = \mu
_1 ^+
(M^* ).
\end{eqnarray*}
By property (\ref{prop2}) this also implies (\ref{1.9}) for sets $M\in
{\mathcal M} ^1 $, proving the first statement of the Theorem.
\\
Now we assume that equality  holds in
(\ref{1.9}) and that  $J$ is convex.
W.l.o.g. we may assume that also $\mu _1 (M) >0$.
For $r>0 $, let $\theta _r (M,x) := \mu _1 (M\cap B_r (x) ) [\mu _1 (B_r
(x))] ^{-1}
$ and  define the upper density of the set $M$ at $x$ by $\theta (M,x):=
\limsup _{ \, r\searrow 0} \theta _r (M,x)$.
Suppose that $M'$ is the set of upper density points of $M$, i.e.
$M'= \{ x\in \R :\, \theta (M,x)=1\} $. Since $\mu _1 (M\Delta M' )=0
$, it
is sufficient to study $M'$ instead of $M$. We first claim that
$M' $ is convex.
\\
Assume that $M'$ is not convex.
Then there are points $x_i \in \R $, $i=1,2,3$, $x_1 <x_3 <x_2 $,
such that
$\theta (M,x_1)=\theta (M,x_2 )=1$ and $\theta (M, x_3 )=0$. Then for
$ r<
(1/4 )
\min \{ (x_3 -x_1); (x_2 -x_3 )\} $, and small enough,
we have that $\theta _r (M, x_i ) \geq
(7/8) $, $i=1,2$, and $\theta _r (M,x_3 )\leq (1/8)$. Let $\{ M_k \} $ a
sequence
of smooth sets such that $\mu_1 (M_k \Delta M) \to 0$ and $P_{\mu
_1 } (M_k
)
\to P_{\mu _1 } (M) $, as $k\to \infty $. For $k$ large enough - say
$k\geq k_0 $ - we still have $\theta _r (M_k ,x_i ) \geq (3/4)$,
$i=1,2$,
and
$\theta _r (M_k , x_3 ) \leq (1/4) $. This implies that the sets $B_r
(x_i
)\cap M_k $, $i=1,2$, and $B_r (x_3 )\setminus M_k $ are nonempty for
these
$k$. In other words, if $k\geq k_0 $, then there is a nonempty interval
$I_k =(y_k ^1 , y_k ^2 ) $ such that $I_k \subset (x_1 -r ,x_2 +r )$,
$I_k
\cap M_k =\emptyset
$ and
$y_k ^1 ,
y_k ^2 \in \partial M_k $. Setting $N_k := I_k \cup M_k $ we then
have in
view of the isoperimetric property (\ref{1.9}),
\begin{eqnarray*}
P_{\mu _1 } (M_k ) - J(\mu _1 (M_k )  ) & \geq &
P_{\mu _1 } (N_k ) - J(\mu _1 (N_k ) ) +
\psi (y_k ^1 ) + \psi (y_k ^2 )
\\
    & \geq   &
\psi (y_k ^1 ) + \psi (y_k ^2 ) \geq \delta , \quad \forall k\geq k_0 ,
\end{eqnarray*}
for some $\delta >0 $ independent on $k$. Passing to the limit for
$k\to \infty $, this also implies $P_{\mu _1 } (M) -J(\mu _1 (M)) \geq
\delta >0$, that is a contradiction. Hence $M' $ is convex.
\\
Now we assume that equality  holds in
(\ref{1.9}) and that  $J$ is strictly convex. Hence $M'=(a,b)$, where $a,b\in \R $, $a<b$.
Setting $\alpha =\Psi (a)$, $\beta = \Psi (b) $, we have that
$M^* =(-c,c)$,
where $\Psi (c) = (\beta -\alpha ) /2$, so that
$$
J(\alpha ) +J(\beta ) = 2J\Big( \frac{\beta -\alpha }{2} \Big) .
$$
The strict convexity of $J$ then implies that $|\alpha |= |\beta |$,
that
is $b=-a$. The Theorem is proved.
$\hfill \Box $

\section {Product measures on $\R^n$}
In this section we prove isoperimetric inequalities
with respect to product measures on $\R^n$
and we apply them to obtain integral inequalities in Sobolev spaces.

\noindent We deal with a product measure $\mu$ on $\R ^n $ defined by
\begin{equation}\label{nuovo}
d\mu =\psi (x_1 )\rho (x' ) \, dx \,,
\end{equation}
where  $x=(x_1 ,x_2 ,\ldots ,x_n )= (x_1 ,x')$
is a point in $\R^n $, ($n\geq 2$),
$\psi $ is a function as in the previous section
and $\rho$ is a positive, continuous function on $\R^{n-1}$.

\subsection {Isoperimetric inequalities}

If $M\subset \R ^n $
we introduce $x'$- slices, $x'\in \R^{n-1 } $, by
$$
M (x') = \{ x_1 : \, (x_1 ,x') \in M\}.
$$
Notice that if $M\in {\mathcal M} ^n $, and if $\mu (M)<\infty $
then $\mu _1 (M(x') )<\infty  $ for almost
every $x' \in \R ^{n-1}$.
For  $M\in {\mathcal M} ^n $,  we define its {\sl Steiner
$\mu $-symmetrization\/} $M^* $ by
\begin{equation}
\label{1.6000}
M^* := \Big\{ x=(x_1 ,x' ) : \,
x_1 \in \Big( M(x') \Big) ^* ,\ x' \in \R ^{n-1}  \Big\} .
\end{equation}
Then
\begin{equation}
\label{1.300}
\mu (M)=\mu (M^* ),
\end{equation}
and it is easy to see that the monotonicity properties
(\ref{1.101})-(\ref{1.103})
carry over to Steiner $\mu $- symmetrization.

As in the one-dimensional case, the above definitions will be read
pointwise
for open and closed sets.
\\
If $M$ is open (respectively closed) then the function
$\varphi (x') := \mu (M(x'))$, ($x'\in \R ^{n-1} $),
is lower (respectively upper) semicontinuous. Since the sets
$(M(x' ))^* $, ($x' \in \R ^{n-1} $),
are open (respectively closed) it then follows that $M^* $
is open (respectively closed), too.
\\
Steiner-like rearrangements
are characterized by the fact that the isoperimetric property (\ref
{1.13})
{\sl on slices\/} carries over to sets on $\R ^n $ (see \cite{[BuZa]}). By adapting
the proof of Theorem \ref{th1} in \cite{[BBMP2]}, we can prove the following result
\begin{mylem}\label{lem1}
The property (\ref{1.13}) holds for Borel sets $M$ on $\R ^n $.
\end{mylem}
By the definition of $\mu ^+ $ and $P_{\mu} $, and by property (\ref
{prop2})
we then also have

\begin{mycor}
The following inequalities hold
\begin{eqnarray}
\label{isopMink}
\mu ^+ (M) & \geq & \mu ^+ (M^* ) \quad \mbox{for every Borel set
$M$, and }
\\
\label{isoperi}
P_{\mu } (M) & \geq &
P_{\mu } (M ^* ) \quad \mbox{ for every $M\in {\mathcal M}
^n
$.}
\end{eqnarray}
\end{mycor}

The next property follows easily from Lemma \ref {lem1}.

\begin{mycor}\label{cor2}
Let $M,N$ opens sets in $\R ^n $ with $M\subset N$. Then
\begin{equation}
\label{1.400}
\mbox{dist } \{ M; \partial N\} \leq
\mbox{dist } \{ M^* ; \partial N^* \} .
\end{equation}
\end{mycor}
\bigskip

Finally we analyze the equality case in (\ref{isoperi}).
The following result holds

\begin{myth} \label{th2}
Assume that equality holds in (\ref{isoperi}) for some $M\in
{\mathcal M}
^n $. Then $M(x')$ is either empty or equivalent to an interval for
almost
every
$x'\in
\R ^{n-1}$. Moreover, if
$J$ is strictly convex, then $M=M^* $.
\end{myth}

The proof of Theorem \ref{th2} depends on a precise estimate for the deficit
of the perimeter under Steiner symmetrization for polyhedra. This
approach
is well-known in the case of  the uniform Lebesgue measure $\varphi
\equiv
1$ (see \cite{[BuZa]}, chapter 14).
Let us first introduce
some notation. Let $\nabla ' $ denote the vector of derivatives
$(\partial /\partial x_2 ,\ldots ,\partial /\partial x_n )$.
If $x'\in
\R ^{n-1}$ then let
$l_{x'} $ denote the line $\{ (t,x'):\, t\in \R \} $.
    Let
${\mathcal P }$ denote  the set of polyhedra $\Pi $  in
$\R ^n
$ such that
$l_{x' } \cap \Pi $ is either empty or consists of a finite number of
points for every $x' \in \R ^{n-1}$. The map $p: \partial \Pi \to \R ^
{n-1}$
will be called a projection. If $\Pi \in {\mathcal P} $ then $\R ^{n-1} $ is
splitted into a finite numbers  of domains $Q$ such that the part of
$\partial \Pi $ which is projected into $Q$ consists of a finite
number $2m
$ of components $\Gamma _j $ whose projections onto $Q$ are one-to-
one. (The
number $m$ depends on $Q$, and those for which $m=0$ will not be
considered
further on.)  Each  $\Gamma _j $ permits an (affine) representation
$x_1 = z_j (x') $, $x'\in Q$. Then it follows that
\begin{eqnarray}
\label{polyh1}
P_{\mu }(\Pi ) & = & \sum\, ^* \int\limits_Q \sum\limits_{j=1} ^{2m}
\sqrt{
1+|\nabla ' z_j | ^2 } \, \psi (z_j ) \rho (x') \, dx',
\end{eqnarray}
where the sum $\sum ^*$ is taken over all the $Q$ for which $m\geq 1 $.
After the Steiner symmetrization, the area of the boundary will be
\begin{eqnarray}
\label{polyh2}
P_{\mu }(\Pi ^* ) & = & \sum\, ^* \int\limits_Q 2 \sum\limits_{j=1} ^
{2m}
\sqrt{ 1+|\nabla ' z | ^2 } \, \psi (z ) \rho (x') \, dx',
\end{eqnarray}
where the function $z: Q\to \R_0 ^+ $ is given by
\begin{equation}
\label{polyh3}
\sum\limits_{j=1}^{2m} (-1) ^j \Psi (z_j ) = 2\Psi (z).
\end{equation}
\begin{mylem}
Let $\Pi \in {\mathcal P}$.
Then, with the above notations,
\begin{eqnarray}
\nonumber
    & & P_{\mu } (\Pi ) -P_{\mu } (\Pi ^*)
\\
\label{polyh4}
& \geq &
\Big( P_{\mu } (\Pi ^* ) \Big) ^{-1} \sum \, ^* \int\limits_{Q} \sqrt
{\psi
(z) \Big| \sum_{j=1} ^{2m} \psi (z_j ) -2\psi (z)\Big| } \rho (x') \,
dx' .
\end{eqnarray}
\end{mylem}
{\bf Proof:\/}
For convenience,
we set $y_j := \Psi (z_j ) $, $j=1,\ldots 2m$, $y=\Psi (z)$, and $J:=
\psi
(\Psi ^{-1})$. Recall that $J$ is convex by our assumptions. Then we
find
\begin{eqnarray}
\nonumber
    & & \sum\limits_{j=1} ^{2m} \sqrt{
1+|\nabla ' z_j | ^2 } \, \psi (z_j ) -2 \sqrt{1+ |\nabla ' z|^2 }
\psi (z)
\\
\nonumber
    & = & \sum\limits_{j=1} ^{2m} \sqrt{ J(y_j ) ^2 + |\nabla ' y_j |
^2 } -
2\sqrt{J(y) ^2 +|\nabla ' y| ^2}
\\
\nonumber
    & \geq & \frac{J(y)}{\sqrt{J(y)^2 +|\nabla ' y|^2 }}  \Big(
\sum\limits_{j=1} ^{2m} J(y_j ) - 2J(y)\Big)
\\
\label{polyh5}
    & = &
     \frac{ \sum_{j=1} ^{2m} \psi (z_j ) -2\psi (z)}{\sqrt{ 1
+|\nabla ' z|^2 }}  .
\end{eqnarray}
Integrating (\ref{polyh5}) and applying Cauchy-Schwarz inequality
we have (\ref{polyh4}).
$\hfill \Box $
\\[0.3cm]
{\bf Proof of Theorem \ref{th2}:\/} For $x\in \R ^n $ and $r>0$ let $\theta _r
(M,x) := \mu (M\cap B_r (x)) [\mu (B_r (x))] ^{-1} $, and define $\theta
(M,x') $ and the set $M' $ of upper density points of $M$ as in the
proof of Theorem 1. As before, we may restrict ourselves to the set
$M'$  instead of $M$.
\\
Choose a sequence of polyhedra $\{ \Pi _k \} $ such that
$\mu (\Pi _k \Delta M)\to 0$ and $P_{\mu } (\Pi _k ) \to P_{\mu } (M)
$ as
$k\to
\infty $. Without loss of generality, we may also assume that $\Pi _k \in {\mathcal P}$,
$k=1,2,\ldots $. Since $\mu (\Pi _k ^* \Delta M^* ) \to 0$ as $k\to
\infty$, we have by the lower semicontinuity of the perimeter,
\begin{equation}
\label{polyh6}
\lim_{k\to \infty } (P_{\mu } (\Pi _k ) -P_{\mu } (\Pi _k ^* )) =0.
\end{equation}
Set
$$
R_k := \{ x' \in \R ^{n-1} :\, l_{x'} \cap \Pi _k  \mbox{ has at least
two components }\},
$$
and introduce a measure $\nu $ on $\R ^{n-1} $ by
$$
d\nu = \rho (x')\, dx' .
$$
Since the function $\psi $ is bounded away from $0$, the previous Lemma
2 together with (\ref{polyh6}) tells us that
\begin{equation}
\label{polyh7}
\lim_{k\to \infty } \nu (R_k ) =0 .
\end{equation}
We claim  that $M(x')$ is convex for almost every $x'\in
\R ^{n-1} $. \\
Assume that this is not the case. Then there are points
$x^i = (z_i , x'_0 )$, $i=1,2,3$, with $x' _0 \in \R ^{n-1} $, and
$z_1 <z_3
<z_2
$, such that $\theta (M,x^1 )=\theta (M,x^2 )=1$, and
$\theta (M, x^3 )=0$.
Let $\varepsilon >0$ and small (the exact choice of $\varepsilon $
being specified later).
Choose $r (=r(\varepsilon ))>0 $, and small enough such
that $\theta _r (M, x^i ) \geq 1-\varepsilon $, $i=1,2$, and $\theta
_r (M,x^3 )
\leq \varepsilon$.
For $k$ large enough - say $k\geq k_{\varepsilon}$ - we then still have
$\theta _r (\Pi _k ,x^i ) \geq 1-2\varepsilon $, $i=1,2$, and $\theta
_r (\Pi _k ,x^3
) \leq 2 \varepsilon $.
Let
$$
H_k := \{ x=(x_1 ,x' ) \in B_r (x^3 ) :\, \nu (l_{x'} \cap \Pi _k
\cap B_r (x ^1 )) >0 , \mbox{ and }  \nu (l_{x'} \cap \Pi _k \cap B_r
(x ^2 )) >0 \}.
$$
By choosing $\varepsilon $ small enough we can achieve that $\mu (H_k
) >(1/2 ) \mu ( B_r (x^3 )) $, and in view of $\theta _r (\Pi _k ,
x^3 ) \leq 2\varepsilon $, also that
$\mu (H_k \setminus \Pi _k ) >(1/4) \mu (B_r (x^3 ))$. Hence there is
a number $c_0 >0$ which depends only on $\varepsilon $, but not on
$k$, such that $\nu (R_k \geq c_0 $. But this contradicts
(\ref{polyh7}).
    This proves the claim.
\\
Hence there is a nullset $N\subset \R ^n$,  a measurable set $G
\subset \R
^{n-1} $, and measurable functions
$z_i$,
$i=1,2$, such that
$$
M= N\cup \{(x_ 1 ,x') :\, z_1 (x' )< x_1 <z_2 (x'),\, x' \in G\} .
$$
Using Lemma 2 and the limit property (\ref{polyh6})
 we have
\begin{eqnarray}
\nonumber
    & &
0 =P_{\mu } (M) -P_{\mu } (M^*)
\\
\label{polyh9}
    & \geq &
\Big( P_{\mu } (M)\Big) ^{-1}
\int\limits_G
\sqrt{ \psi (z) \Big| \sum_{j=1} ^{2} \psi (z_j ) -2\psi (z)\Big| }
\rho (x')
\, dx' ,
\end{eqnarray}
where $z$ is given by $2\Psi (z) = \Psi (z_2 ) -\Psi (z_1 )$.
    Using the strict convexity of $J$ this implies that $z_2 = -z_1 =z
$ on
$G$ and the Theorem is proved.
$\hfill \Box $
\subsection{Steiner  $\mu$-symmetrization of functions}

If $u\in {\mathcal F} ^n $ we define its \\
{\sl Steiner}  $\mu$-{\sl symmetrization} {\sl (w.r.t. $x_1 $)\/} $u^* $ by
\begin{equation}
\label{2.1u*}
u^* (x) := \sup \Big\{ t\in \R :\, x\in \{ u>t \} ^* \Big\} ,\quad x
\in \R ^n .
\end{equation}
By its definition, the function $u^* $ is nonincreasing and
right-continuous
     w.r.t.
the variable $x_1 $. Moreover $u$ and $u^* $ are {\sl equimeasurable\/} functions,
that is
\begin{equation}
\label{2.2u*}
\{ u>t \} ^* = \{ u^* >t \} \ \mbox{ and } \
\{ u\geq t \} ^* = \{ u^* \geq t \} \quad \forall t>\inf u ,
\end{equation}
which implies that $\forall t>\inf u $ and for a.e. $x' \in \R ^{n-1} $,
\begin{eqnarray}
\nonumber
     & & \mu \Big( \{ u>t \} \Big) =
\mu \Big( \{ u^* >t \} \Big) ,\
\mu \Big( \{ u\geq t \} \Big) =
\mu \Big( \{ u^* \geq t \} \Big) ,\ \mbox{ and }\\
\nonumber
     & & \mu \Big( \{ u(\cdot ,x') >t \} \Big) =
\mu \Big( \{ u^* (\cdot ,x' ) >t \} \Big) ,\\
\label{2.3u*}
     & &
\mu \Big( \{ u (\cdot ,x' )\geq t \} \Big) =
\mu \Big( \{ u^* (\cdot ,x' )\geq t \} \Big).
\end{eqnarray}
 Furthermore,
the monotonicity
(\ref{1.101}) implies
\begin{equation}
\label{2.31}
u,v\in {\mathcal F} ^n ,\ u\leq v\Longrightarrow u^* \leq v^* .
\end{equation}
\medskip
\begin{myrem} \label{rem5}\rm
We will generally
not distinguish between $u$ and its equivalence class given by
all measurable functions which differ from $u$ on a nullset. But if
$u$ is continuous, then the sets $\{ u>t\} $,
($t\in \R $) are open,
and we will  agree that
the above definition of $u^*$ has to be understood in pointwise sense.
Furthermore it is easy to see that the
sets $\{ u^* >t\} $,
(respectively $\{ u^* \geq t\} $) are open (respectively closed),
($t\in \R $),
so that $u^* $ is continuous too.
\end{myrem}
\medskip

\begin{myrem} \rm An equivalent
definition of $u^* $ can be given by using
the $\mu $- distribution function of $u$
(w.r.t.
$x_1 $),  $m_u $,  defined by
$$
m_u (x' ,t) := \mu _1 \Big( \{ u(\cdot ,x' )>t \} \Big) ,
\quad x' \in \R ^{n-1} ,\ t>\inf u .
$$
The function $m_u $ is nonnegative,  nonincreasing,  right-continuous
w.r.t. the variable $t$, and
\begin{equation}
\label{2.4u*}
u^* (x_1 ,x') = \sup \{ t\in \R :\, m_u (x' ,t)> \mu _1 ((-
x_1 ,x_1 )) \} ,
\quad
x=(x_1 ,x' )\in \R ^n .
\end{equation}
\end{myrem}
\hspace*{1cm}

Proceeding analogously as we did in \cite{[BBMP2]} for a
certain class
of Steiner-like rearrangements w.r.t. a finite measure, we can prove
the
following properties.
\begin{myth} \label {th3}
$\mbox{ }$
\\
{\bf 1)}
If
$u\in L^1 _+ (\mu , \R ^N)$, then
\begin{equation}
\label{2.8}
\int\limits_{\R ^n } u\, d\mu = \int\limits_{\R ^n} u ^* \, d\mu .
\end{equation}
{\bf 2)} If $u\in {\mathcal F} ^n$ and if
$\varphi :\, \R \rightarrow \R $ is a nondecreasing function, then
\begin{equation}
\label{2.9}
     \varphi (u^* )
= \Big( \varphi (u)\Big) ^*.
\end{equation}
{\bf 3)} (Cavalieri's principle)
\\
If $f:\R \rightarrow \R$ is continuous or nondecreasing,
$u\in {\mathcal F} ^n $ and if
$f(u)\in L^1 (\mu ,\R ^n)$,  then
$f(u^* )\in L^1 (\mu ,\R ^n )$ and
\begin{equation}
\label{2.11}
\int\limits_{\R^n } f(u)\ d\mu  =\int\limits_{\R^n } f(u^* )\ d\mu .
\end{equation}
{\bf 4)}
Let $F\in C((\R _0 ^+ )^2 )$, $F(0,0)=0$, and
\begin{equation}\label{2.14}
F(A,B)-F(a,B) -F(A,b)+F(a,b)   \geq   0
\end{equation}
for all  $a,b,A,B\in \R
_0 ^+$
with $a\leq A,\  b\leq B.$

\noindent Furthermore, let $u,v\in {\mathcal F} ^n _+ $ such that $F(u,0)$,
$F(0,v)$,
$F(u,v)\in L^1 (\mu ,R ^n )$.
Then
\begin{equation}
\label{2.15}
\int\limits_{\R ^n } F(u,v)\, d\mu  \leq \int\limits_{\R ^n } F(u^* ,v^*
)\, d\mu .
\end{equation}
{\bf 5)}
(Nonexpansivity of the rearrangement )
\\
Let $G\in C(\R _0 ^+ )$ continuous, nondecreasing and convex
with $G(0)=0$,
and let $u,v\in {\mathcal F} ^n _+$ such
that
$G(|u|),G(|v|), G(|u-v|) \in L^1 (\mu ,\R
^n )$. Then
\begin{equation}
\label{2.21}
\int\limits_{\R ^n } G(| u^*-v^* |)\,
d\mu \leq \int\limits _{\R ^n } G(|u-v|) \,
d\mu.
\end{equation}
{\bf 6)}
Let $u,v\in L^2 _+ (\mu ,\R ^n )$.
Then
\begin{equation}
\label{2.22}
\int\limits_{\R^n} uv\, d\mu \leq \int\limits_{\R^n } u^* v^* \, d\mu.
\end{equation}
{\bf 7)}
Let $u,v\in  L^{\infty } (\R ^n )\cap {\mathcal F} ^n $.
Then
\begin{equation}
\label{2.23}
\Vert u^* -v^*\Vert _{\infty } \leq \Vert u-v\Vert _{\infty } .
\end{equation}
\end{myth}
\begin{myrem}\label{remark4} If $F\in C^2 $ then (\ref{2.14}) is equivalent to
\begin{equation}
\label{2.16}
\frac{\partial ^2 F (\sigma ,\tau ) }{\partial \sigma  \partial \tau }
\geq 0,
\quad
\forall \sigma ,\tau \in \R _0 ^+ .
\end{equation}
\end{myrem}
\medskip

\noindent {\bf Proof of Theorem \ref{th3}:\/}
     The proofs of the properties {\bf 1)}-{\bf 3)}, and of
{\bf 5)}-{\bf 7)}
     mimic
the proofs of analogous properties in \cite{[BBMP2]}.
\\
{\bf 4)}  We proceed similarly as in \cite{[CrRoZw]}.
In view of (\ref{2.14}),
there exists a nonnegative measure,  denoted by
$dF_{\sigma \tau }$, such that
\begin{eqnarray}
\nonumber
     & & F(s,t)-F(s,0)-F(0,t) \  = \
\int\limits_0 ^t \int\limits _0 ^s
\, dF_{\sigma \tau } \\
     & = & \int\! \! \! \int\limits_{(\R _0 ^+ )^2 }
\chi (0,s) (\sigma  ) \, \chi (0,t) (\tau )
     dF_{\sigma \tau }
\label{2.17}
\qquad
\forall s,t
\in
\R _0 ^+ .
\end{eqnarray}
Notice that in the case $F\in C^2 $ we have
$$
dF_{\sigma \tau } = \frac{\partial ^2 F (\sigma  ,
\tau )}{\partial \sigma
\partial \tau }  \, d\sigma  \, d\tau  ,
$$
which is obviously nonnegative by (\ref{2.14}).\\
Choosing $s=u(x)$ and $t=v(x)$ in (\ref{2.17}) and then
integrating we find
\begin{eqnarray}
\nonumber
\int\limits_{\R ^n} F(u,v)\, d\mu   & = &
\int\limits_{\R ^n } F(u,0)\, d\mu +
\int\limits_{\R ^n } F(0,v)\, d\mu
\\
\label{2.19}
     & &
+ \int \! \! \! \! \! \! \! \int\limits_{(\R _0 ^+ ) ^2 }
\mu \Big( \{ u>\sigma  \}\cap \{v>\tau \}\Big)
\, dF_{\sigma \tau }.
\end{eqnarray}
An analogous expression for $\int_{\R ^n } F(u^*,v^*)\, d\mu  $ holds.
Since from (\ref{1.102}), we have
\begin{equation}
\label{2.10000}
\quad \mu \Big( \{ u>\sigma \}\cap \{v>\tau \}\Big)
 \leq  \mu \Big( \{ u^*>\sigma \}\cap \{v^*>\tau \}\Big) \qquad
\forall \sigma ,\tau >0.
\end{equation}
Then (\ref{2.15}) follows from (\ref {2.19}) and (\ref {2.10000}).
$ \hfill \Box $
\medskip

\begin{myrem}\rm
Let $M=\R \times M' $, where $M\in {\mathcal M} ^{n-1} $,
and let
$u,v\in {\mathcal F} ^n $. Then the properties {\bf 1)}, and {\bf 3)}-
{\bf 6)}
of Theorem \ref{th3}
hold with the range of integration restricted to $M$.
Indeed the proof of this result can be easily obtained by
Theorem \ref{th3} setting
$$
U(x):= \left\{
\begin{array}{ll}
u(x) & |\ \mbox{ if } \ x\in M\\
\inf u & | \ \mbox{ if } \ x\in \R ^n \setminus M
\end{array}
\right.
$$
and therefore
$$
U^* (x):= \left\{
\begin{array}{ll}
u^* (x) & |\ \mbox{ if } \ x\in M\\
\inf u & | \ \mbox{ if } \ x\in \R ^n \setminus M.
\end{array}
\right.
$$
\end{myrem}
\medskip

We conclude this subsection with the following property, which  is easy
to prove but it is crucial
for the next section.

Let $u$ be a function belonging to $ C(\R ^n )$. Denote the
{\sl modulus of continuity of $u$\/}, by
\begin{equation}
\label{2.5}
\omega _u (t) := \sup \{ |u(x)-u(y)|:\, |x-y|<t\} , \quad t>0.
\end{equation}
Notice that $u$ is uniformly continuous iff
$ \, \lim_{t\searrow 0} \, \omega _u(t) =0$.

\begin{myprop}
Let $u\in C(\R ^n )\cap {\mathcal F} ^n $. Then $u^* \in C(\R ^n )$ and
\begin{equation}
\label{4.1}
\omega _u (t)\geq
\omega _{u^* } (t) \quad \forall t>0.
\end{equation}
In particular, if $u$ is Lipschitz continuous
with constant $L$, then also   $u^* $ is Lipschitz continuous
with a Lipschitz constant  $L^*$ such that $L^* \le L$.
\end{myprop}

\noindent{\bf Proof : \/}   By Remark \ref{rem5},  and
by Corollary \ref{cor2}, (\ref{1.400}),  since $u$ is continuous we have that
\begin{equation}
\label{4.2}
\mbox{dist } \Big\{ \{u>t\} ;\partial \{ u>s\} \Big\} \leq
\mbox{dist } \Big\{ \{u^* >t\} ;\partial \{ u^* >s\} \Big\}\quad
\forall s,t\in \R
\mbox{ with $s<t$,}
\end{equation}
which implies (\ref{4.1}).
Since $L= \sup \{ \omega _u (t)/t :\, t>0 \} $, and similarly for
$u^* $,
the second assertion follows, too.
$\hfill \Box $

\subsection{Integral inequalities in Sobolev spaces}

In this subsection we state  integral inequalities which
involve derivatives of a function and its rearrangement.
Variants of them are well-known for different types of rearrangements, including Steiner
symmetrization (see, for instance, \cite{[Ka],[CiFu],[ET]}),
and they are usually referred as Polya-Szeg\"o type inequalities.
Theorem \ref{ps} below can be shown as Theorem 6.1 and Corollary 6.1 in
\cite{[BBMP2]}, and its proof is therefore omitted.
\begin{myth} {\sl (Polya-Szeg\"o's principle)\/}\label{ps}
\\
Let $G=G(y,v,x') $ be a function belonging to
$L^{\infty } ( \R ^n \times \R_0 ^+ \times \R^{n-1} ) $ where $y=(y_1
,\ldots ,y_n ) \in \R ^n $). Assume also that
$G $ is continuous in $v$, convex in $y$, even in
$y_1 $ and nondecreasing in $y_1$ with $y_1>0$. Moreover  let
$u\in L^{\infty } (\R ^n ) $ a nonnegative
Lipschitz continuous function with compact support. Then
\begin{equation}
\int\limits_{\R^n } G( \nabla u ,u,x' ) \, d\mu \, \geq \,
\int\limits_{\R^n } G( \nabla u^* ,u^* ,x')\, d\mu  .
\label{5.4.11}
\end{equation}
Moreover, if $u\in W ^{1,p} _+ (\mu ,\R ^n )$, for some $p\in [1,
\infty )$,
then $u^* \in
W^{1,p} _+ (\mu ,\R ^n )$, and inequality (\ref{5.4.11}) holds if
\begin{equation}
\label{5.4.110}
|G(y,v,x')|\leq C|y|^p \quad
\mbox{ for some $C>0$,}
\end{equation}
    for any  $ (y,v,x')\in \R ^n \times
\R
\times \R ^{n-1} $.
Finally, if $\Omega $ is a bounded domain in
$\R ^n $ and $u$ is a nonnegative function belonging to $W_{0} ^{1,p}
(\mu ,\Omega )$
for some $p\in [1,\infty )$, then
$u^* \in W^{1,p} _{0} (\mu , \Omega ^*)$.
\end{myth}
Using Theorem \ref{th2} and proceeding analogously as in \cite{[BBMP2]},
proof of Theorem 5.4, one can obtain a criterion for the equality
case in
the inequality (\ref{5.4.11}). We omit the proof.
\begin{myth}\label{3.5}
Let $u\in W^{1,p} (\varphi ,\R ^n ) $ for some $p\in [1,\infty )$.
Furthermore, let the function $J := \psi (\Psi ^{-1})$ be strictly
convex, let $G \in C(\R ^n )$, $G=G(y)$, $y=(y_1 , \ldots ,y_n )$,
$G$ convex and strictly increasing in $y_1 $ for $y_1 >0$, and such
that
\begin{equation}
\label{5.5.100}
|G(y)|\leq C(1+|y| ^p ), \quad \mbox{for some $C>0$.}
\end{equation}
Finally, assume that
\begin{equation}
\label{5.5.2}
\int\limits_{\R ^n } G(\nabla u)\, d\mu =
\int\limits_{\R ^n } G(\nabla u^*)\, d\mu .
\end{equation}
Then $u=u^{*} $.
\end{myth}

\section{Radial measures}

In this section we consider  measures $\mu $ whose density is a radially
symmetric function i.e.
\begin{equation}
\label{4.800}
d\mu = \varphi (|x|) dx,
\end{equation}
where $\varphi \in C(\R _0 ^+ )$ is positive.

\noindent We prove  isoperimetric inequalities with respect
two special measures whose densities are
$$
\varphi_1(|x|)=\mbox{exp}\, \{ c|x|^2 \},\quad c>0,
,\qquad
 \varphi_2(|x|)= |x|^{1-n} \mbox{exp } \{ a(|x|) \},\qquad
$$
where $a\in C(\R _0 ^+ )$ is convex.
\\

If $M\in {\mathcal M} ^n $, and if $\mu (M) <+\infty $, then
let $M^{\displaystyle \star} $ denote the ball $B_R $ such
that  $\mu (M) =\mu (B_R )$. If $\mu (M)=+\infty $ then let
$M^{\displaystyle
\star } =\R ^n $. We call $M^{\displaystyle \star }$ the
{\sl Schwarz $\mu$-symmetrization\/} of $M$. As in the previous
sections,
we replace this
definition by
pointwise ones for open and closed sets.
Thus, if $M$ is open/closed with finite
$\mu$-measure, then let $M^{\displaystyle \star } $
the open/closed ball centered
at zero, having the same measure as $M$.
\\
We  ask for additional conditions on the measure $\mu $ such that
the following isoperimetric inequality holds
\begin{equation}
\label{4.1000}
P_{\mu } (M) \geq P_{\mu } (M^{\displaystyle \star }), \quad \forall
M\in {\mathcal M} ^n ,
\end{equation}
with equality iff $M= M^{\displaystyle \star} $.
Although we are not able to give a necessary and sufficient condition
for (\ref{4.1000}) to hold, we show below that the above isoperimetric
property holds if
\begin{equation}\label{4.1999}
d\mu =\mbox{exp}\, \{ c|x|^2 \}dx,
\quad c>0
\end{equation}
or if
\begin{equation}
\label{5.1000}
d\mu = |x|^{1-n} \mbox{exp } \{ a(|x|) \} \, dx, \quad
\mbox{where $a\in C(\R _0 ^+ )$ is convex.}
\end{equation}

\begin{myth}  \label{th6}
Let $\mu $ be the measure defined by (\ref{4.800}) with $\varphi$ defined by \eqref{4.1999}.
Then
\begin{equation}
\label{4.20}
M^{\displaystyle \star } + B_r  \subset
\Big( M+ B_r  \Big) ^{\displaystyle \star },
\quad \mbox{$\forall $ Borel sets and $\forall r>0$,}
\end{equation}
and
(\ref{4.1000}) holds. Furthermore, if $M,N$ are open sets with $M
\subset N$,
then
\begin{equation}
\label{4.21}
\mbox{dist }\{ M; \partial N\} \leq \mbox{ dist }
\{ M^{\displaystyle \star} ; \partial N^{\displaystyle \star } \} .
\end{equation}
Finally, if $P_{\mu } (M)= P_{\mu }
(M^{\displaystyle \star }) $ for some $M\in {\mathcal M} ^n $, then $M=
M^{\displaystyle \star } $.
\end{myth}
\noindent{\bf Proof :\/} Let $M$ compact, $M\subset \overline{B_R }$ for some
$R>0$, and set
\begin{eqnarray*}
A(M) & := &
\{ N\subset \R ^n :\, N \mbox{ compact }, \, N\subset
\overline{B_R } ,\, \mu (N) =\mu (M),
\\
     & & \qquad \mu (N +B_r )
\leq
\mu  (M+ B_r ) \ \ \forall r>0 \}.
\end{eqnarray*}
Letting
$$
\delta := \inf \{ \mu (N\Delta M^{\displaystyle \star } ):\,
N\in A(M)\} ,
$$
there exists a sequence $\{ N_k\} \subset A(M) $ with
$$
\lim_{k\to \infty } \mu (N_k \Delta M^{\displaystyle \star } ) =
\delta .
$$
Since the $N_k $'s are equibounded, there is a subsequence
$\{ N_{k' } \} $ which converges in Haussdorf distance to some
set $N$,
which also implies that $N\in A(M) $ and
$\mu (N\Delta M^{\displaystyle \star } )=\delta $.
Assume that $\delta >0$. Then we find two density points $\eta ,\zeta $
of $N$ and
$M^{\displaystyle \star }$ such that $\eta \in M^{\displaystyle \star }
\setminus N$ and $\zeta \in N\setminus M^{\displaystyle \star } $.
After some
rotation of the coordinate system
$$
x=\rho (\xi ), \quad
(x,\xi \in \R ^n ),
$$
we achieve
$$
\rho (\eta )= y=(y_1 ,y' ), \ \rho (\zeta ) =z =(z_1 ,y' ),
$$
for some $y'\in \R ^{n-1} ,y_1 ,z_1 \in \R $. Defining $N' $ by
$$
N' :=
\rho (N) \equiv  \{ x=\rho (\xi ):\, \xi \in N\}
$$
let $(N')^*$ denote its Steiner $\mu $-symmetrization w.r.t. the
variable $x_1 $. Notice that $\mu ( N' \Delta M^{\displaystyle \star
} ) =\delta $, and,  due to the product structure of $\varphi $, we
have that  $(N')^{\displaystyle \star } = ((N')^* )^{\displaystyle
\star } =N^{\displaystyle \star}$, and $N', (N')^* \in A(M)$. Since
the slices $(N'(y' ))^* $ and $M ^{\displaystyle \star } (y') $ are
intervals centered at zero, it is easy to see that
$$
\mu _1 \Big(
(N' (y' ) )^* \cap M^{\displaystyle \star }
(y') \Big) >
\mu _1 \Big(
N' (y' )  \cap M^{\displaystyle \star }
(y')  \Big) .
$$
This also implies
$$
\mu \Big(
(N' )^* \cap M^{\displaystyle \star }
\Big) >
\mu  \Big(
N'  \cap M^{\displaystyle \star }
\Big) ,
$$
and thus
$$
\mu \Big( (N')^* \Delta M^{\displaystyle \star } \Big)
<
\mu \Big( N' \Delta M^{\displaystyle \star } \Big),
$$
contradicting
the minimality of $\delta $. Hence $N= M^{\displaystyle \star } $,
and (\ref{4.1000}) is proved for compact sets. It is easy to see
that this also implies property (\ref{4.20}) for compact sets,
and by a simple
approximation argument, as well for Borel sets.
It is well-known that (\ref{4.20}) also implies
(\ref{4.21}), and the lower semicontinuity of the perimeter $P_{\mu } $
yields (\ref{4.1000}).
\\
Assume finally that $P_{\mu } (M) = P_{\mu } (M^{\displaystyle
\star } )$
for some $M\in {\mathcal M} ^n $. Let $\rho $ any rotation about the
origin, and let $^* $ denote $\mu $- Steiner symmetrization in
direction $x_1
$. Then $[(\rho M )^* ] ^{\displaystyle \star } = M^{\displaystyle
\star }$,
which means that
$P_{\mu } (\rho M) = P_{\mu } ((\rho M) ^*) $. By Theorem \ref{th2} this
implies that
$\rho M= (\rho M) ^* $. Since $\rho $ was arbitrary
we just have proved that $M$ is symmetric w.r.t. every $(n-1)$-
hyperplane
through the origin. Hence $M= M^{\displaystyle \star} $.
$\hfill \Box$
\\[0.3cm]
\hspace*{1cm}

Next we rewrite the isoperimetric inequality in terms of  $\mu (M)$.  Let
$I_{\mu}(m)$ be the isoperimetric function defined in \eqref{im},
\begin{eqnarray}
\label{4.1100}
h(r) & := & n\omega _n e^{cr^2} r^{n-1} \quad \mbox{and }\\
\label{4.1101}
    H(r) & := & \int_0 ^{r} h(t)\, dt .
\end{eqnarray}
Then
$$P_{\mu } (M^*) = h\big( H^{-1} (\mu (M))\big)=I(\mu(M^*)) .$$
Therefore (\ref{4.1000}) reads as
follows:
\begin{mycor} \label{isosf}
If $\mu (M)  <+\infty $, then
\begin{equation}
\label{4.1102}
P_{\mu } (M) \geq h\big( H^{-1} (\mu (M))\big)=I(\mu(M^*)).
\end{equation}
\end{mycor}
\bigskip

Now let us define the $\mu$-Schwarz symmetrization  of
functions with respect to  the measure $\mu $ defined in
\eqref{4.800} with $\varphi$ defined by \eqref{4.1999}.
First we introduce a function
$\tilde u: ]0,+\infty[ \rightarrow\R$
defined by
$$
\tilde u (s)=\inf\left \{t\in \R:~~m_u(t)\le s    \right \}.
$$
Notice that $\tilde u $ is a
nonincreasing and right-continuous function. Observe also that $
\tilde u (s)$ is the inverse function of
$m_u(t)$, if $ \tilde u (s)$ is not constant on intervals.
In this case, the following equality holds,
\begin{equation}\label{dermu}
\frac{\partial \tilde u (s)}{\partial s}=\left [ \frac{\partial
m_u (t)}{\partial t}  \right]^{-1},
\end{equation}
where $s=m_u (t)$.
\\
If $u\in {\mathcal F} ^n $ we define the
{\sl Schwarz $\mu $-symmetrization of $u$\/}
by
\begin{equation}
\label{4.0}
u^{\displaystyle \star} (x) := \sup \Big\{
t\in \R :\, x\in \{ u>t \} ^{\displaystyle \star} \Big\} ,
\quad x\in \R ^n .
\end{equation}
Observe
that, by definition of $\tilde u$ and $u^\star$, one has
$$
u^{\displaystyle\star}(x)=\tilde u (H(|x|)),
\quad \mbox{for a.e. }   x\in \Omega ^\star . $$

By definition $u^{\displaystyle \star} $ is radially symmetric
and radially decreasing. Moreover  $u$ and $u^{\displaystyle \star }$
are equimeasurable functions.
If $u$ is
continuous then we will understand this definition in
pointwise sense, so that $u^{\displaystyle \star }$ is continuous, too.
The assertions of Theorem \ref{th3}  hold as well for the Schwarz
$\mu $-symmetrization.

As in case of the Steiner $\mu $- symmetrization,
the isoperimetric property
(\ref{4.21}) implies the following estimate for the modulus of continuity.

 \begin {myprop}
Let $\mu $ be given by (\ref{4.800}) and (\ref{4.1999}).
Then
\begin{equation}
\label{4.4}
\omega _u \geq \omega _{u^{\displaystyle \star } } \quad \forall u\in
C(\R ^n )\cap L^{\infty } (\R ^n ) \cap {\mathcal F } ^n .
\end{equation}
\end{myprop}
\medskip

Isoperimetric property
(\ref {4.1000})  and a result due to Talenti \cite{[Ta4]},
imply  some  Polya-Szeg\"o  type inequalties, that is,
integrals involving gradients of a
nonnegative Lipschitz function having compact support decrease under
weighted Schwarz symmetrization. The fact that the equality case in
these inequalities  occurs only in symmetric situations can be shown
 analogously as in \cite{[BBMP2]}, and by
using Theorem \ref{3.5} above. Using arguments as in the proofs of
Theorem \ref{th3} and Corollary \ref{cor2} this leads to norm inequalities in
$W^{1,p} (\R ^n)$. We omit the proofs.
\begin{myth} \label{Polya}
Let $u$ a nonnegative Lipschitz continuous function on $\R ^n $ with
compact support, and let $G\in C(\R_0 ^+ ) $ nonnegative and convex
with $G(0)=0$. Then
\begin{equation}
\label{5.1}
\int\limits_{\R ^n } G(|\nabla u|) \, d\mu  \geq
\int\limits_{\R ^n } G(|\nabla u^{\displaystyle \star} |) \, d\mu .
\end{equation}
Moreover, if $G$ is strictly convex then (\ref{5.1}) holds with
equality sign
only if $u=u^{\displaystyle \star} $. Furthermore, if $u\in W^{1,p} (\mu
,\R ^n )$ is nonnegative, for some $p\in [1,\infty )$, then
$u^{\displaystyle \star} \in W^{1,p} (\mu , \R ^n )$,
and (\ref{5.1}) holds with
$G(t)=t^p $. Finally, if $\Omega $ is a domain in $\R ^n $ and
$u\in W_{0} ^{1,p} (\mu , \Omega )$ is nonnegative, then $u^{\displaystyle
\star } \in W_{0} ^{1,p} (\mu , \Omega ^{\displaystyle \star } )$.
\end{myth}
The isoperimetric inequality leads to several integral inequalities
which
compare an $L^p$- weighted norm   of the gradient of a function
with an $L^q $-weighted
norm of
the same function when the measure $\mu $ is given by
\eqref{4.800} and \eqref{4.1999}.
This type of results are also proved in a
different way in \cite{[MS]}.
\begin{myth}\label{sobolev}
There are constants $C=C(p,q)>0$ such that for every
$u\in C_0 ^{\infty } (\R ^n )$,
\begin{equation}
\label{4.1103}
\Vert \nabla u\Vert _p
\geq C(p,q) \Vert u\Vert _q ,
\end{equation}
where $q\in [p,np/(n-p)]$ for $p\in [1,n)$,
$q\in [n, +\infty )$ for
$p=n$,
and
$q\in [p,+\infty ]$ for $p\in (n,+\infty )$.
Moreover, there are
constants $C (p)>0$ such that for every
$u\in C_0 ^{\infty } (\R ^n )$,
\begin{equation}
\label{4.1104}
\Vert \nabla u\Vert _p
\geq C(p) \Vert u\Vert _{C ^{0, 1- (n/p)} (\R ^n )}  ,
\end{equation}
if
$p\in (n,+\infty )$.
\end{myth}
\noindent {\bf Proof :\/}
\\
{\bf 1)} Let $h,H$ be given by (\ref{4.1100}) and
(\ref{4.1101}). It is then easy to see that
$$
h(r)^q\geq CH(r) \quad \forall r\in [0,+\infty ), \ \ \mbox{ for some
$C>0$,}
$$
if $q\in [1,n/(n-1)]$. Applying Theorem 2.1.1 of \cite{[Ma]},
this implies (\ref{4.1103}) with $p=1$ and $q\in [1,n/(n-1)]$.
\\
{\bf 2)} Let $p>1$, $t>0$, $q\in [1, n/(n-1)]$, $q<p/(p-1)$, and
$u\in C_0 ^{\infty } (\R ^n )$. Applying {\bf 1)} and H\"older's
inequality
we have that
\begin{eqnarray}
\nonumber
\Vert u|u|^{t-1}\Vert _q  & \leq  &
C(1,q) \, \Vert t |u| ^{t-1} |\nabla u|\Vert
_1
\\
\label{4.1105}
    & \leq & C_0 \Vert \nabla u\Vert _p \, \Vert |u|^{t-1} \Vert _{p'} ,
\end{eqnarray}
for some $C_0 >0$, where $p' = p/(p-1)$. Choosing $t =p'/ (p'-q )$, we
obtain
(\ref{4.1103}) for $p\in (1,n)$ with $q\in [p,np/(n-p)]$, and for
$p\geq n$ with $q\in [n,+\infty ) $.
\\
{\bf 3)}
Since $\varphi (t) = \mbox{exp}\, \{ ct^2 \}\ge 1, $ it results
$$
\Vert \nabla u\Vert _p \geq \Bigg( \int_{\R ^n } |\nabla u | ^p \, dx
\Bigg) ^{1/p} ,
$$
we obtain (\ref{4.1104}) from Morrey's Imbedding Theorem. From this we
also obtain (\ref{4.1103}) for $p\in [n ,+\infty )$ with $q=+\infty $.
$\hfill \Box $
\\[0.3cm]
\hspace*{1cm}
Let $X^{p,q} $ denote the closure of $C_0 ^{\infty } (\R ^n ) $ with
respect to the norm
\begin{equation}
\label{4.1106}
|\Vert u\Vert | _{p,q } := \Vert u\Vert _q  +\Vert \nabla u\Vert _p ,
\quad
p,q \in [ 1,+\infty ).
\end{equation}
   From (\ref{4.1103}) one immediately obtains the following results.
\begin{mycor}
$\mbox { } $\\
{\bf 1)}
Let $p\in [1,+\infty ) $, $q\in [p,np/(n-p)]$ for $p<n$, and $q\in
[p,+\infty )$ for $p\geq n$. Then
\begin{equation}
\label{4.1107}
X^{p,q} = W ^{1,p} (\R ^n ,\varphi ).
\end{equation}
{\bf 2)}
Let $\Omega $ any domain in $\R ^n $ and $p\in [1,+\infty )$. Then
$W_0 ^{1,p } (\Omega ,\varphi ) \subset W^{1,p } (\R ^n ,\varphi )$, and
\begin{equation}
\label{4.1108}
\Vert \nabla u\Vert _{p,\Omega } \geq C(p,p) \Vert u\Vert _{p,
\Omega } \quad
\forall u\in W_0 ^{1,p} (\Omega ,\varphi ).
\end{equation}
\end{mycor}
    To our knowledge, the problem of finding the best constants
in the inequalities (\ref{4.1103})   is still open.
Here we solve such a problem in the
special case $p=q=2$:
\begin{mycor}\label{cor3}
Let $c>0$.
Then there holds
\begin{eqnarray}
\label{4.1201}
    & & \inf \Big\{ \frac{\Vert \nabla u\Vert _2 ^2 }{\Vert u\Vert
_2  ^2} : \,
u\in C_0 ^{\infty } (\R ^n )\Big\}  =  2cn
    =     \frac{\Vert \nabla (e^{-c|x|^2 } )\Vert _2 ^2 }{\Vert
e^{-c|x|^2 } \Vert _2  ^2} .
\end{eqnarray}
\end{mycor}
\noindent {\bf Proof:\/} Consider the following eigenvalue problem for the
harmonic
oscillator,
\begin{equation}
\label{4.1202}
-\Delta v +c ^2 |x|^2 v =\lambda  v \quad \mbox{in } \ \R ^n.
\end{equation}
The spectrum and the eigenfunctions are explicitly known (see \cite
{[Tay]},
p.104 ff.). In particular, the spectrum is given by $\{ \lambda =
\lambda _k
:\,
(2k-2+n)c,\, k=1,2,\ldots\} $, the eigenvalue $\lambda _1 = cn$ is
simple and
a corresponding eigenfunction is
$v_1 = \mbox{exp }\{-c|x|^2 /2\} $. This implies
\begin{eqnarray}
\label{4.1203}
    & & \int_{\R ^n } \Big( |\nabla v| ^2 +c^2 |x|^2 |v|^2 -cn |v|^2
\Big) \, dx
     \geq  0
\quad
    \mbox{for every } \ v\in C_0 ^{\infty } (\R ^n),
\\
\nonumber
    & & \int_{\R ^n } \Big( \Big| \nabla (\mbox{exp}\{-c|x|^2 /2 \} )
\Big| ^2
+c^2 |x|^2 \,\mbox{exp}\{-c|x|^2 \} -cn \, \mbox{exp}\{-c|x|^2 \}
\Big)
\, dx  =  0.
\label{4.1204}
\end{eqnarray}
Now let $u\in C_0 ^{\infty} (\R ^n )$. Setting
$v:= u\,\mbox{exp}\, \{c|x|^2 /2\} $ we then find by partial integration
\begin{eqnarray*}
\Vert \nabla u \Vert _2 ^2 -2cn \Vert u\Vert ^2 _2 & = &
    \int_{\R ^n } \Big( |\nabla v| ^2 +c^2 |x|^2 |v|^2 -cn |v|^2
\Big) \, dx,
\end{eqnarray*}
and the assertion follows from (\ref{4.1203}).
$\hfill \Box $
\\[0.3cm]
\hspace*{1cm} Finally we prove an isoperimetric inequality
w.r.t. a measure $\mu $
that is given by
\begin{equation}
\label{5.1000}
d\mu = |x|^{1-n} \mbox{exp } \{ a(|x|) \} \, dx, \quad
\mbox{where $a\in C(\R _0 ^+ )$ is convex.}
\end{equation}
Notice that the measure $\mu $ above is singular at the origin. This
implies in
particular that the outer Minkowski content of the set $\{ 0\} $ is
positive,
namely
$$
\mu ^+ (\{ 0\} ) =n\omega _n e^{a(0)} , \quad (\omega _n :\
\mbox{ measure of the
$n$-dimensional unit ball).}
$$
\begin{mylem}
Let $\Omega $ be a
smooth open set in $\R ^n $ which contains an open neighborhood of
the origin, and such that
$\mu (\Omega )<\infty $, where the measure $\mu $ is given by (\ref
{5.1000}),
and let
$B_R $ the ball with $\mu (B_R ) =\mu (\Omega )$, ($R>0 $). Then
\begin{eqnarray}
\nonumber \mu ^+ (\Omega ) & = &
\int\limits_{\partial \Omega } |x|^{1-n} \mbox{exp } \{ a(|x|)  \} \,
{\mathcal H} ^{n-1} (dx)
\\
\label{5.2}
     & \geq &
n\omega _n e^{ a(R) } = \mu ^+ (B_R ).
\end{eqnarray}
\end{mylem}
\noindent {\bf Proof :\/}
Consider
$$
f(x):= e^{-a(|x|)} \mbox{ div }
\Big( |x|^{-n} e ^{ a(|x|) }
x \Big) =
$$
and denote by $\nu $ the  exterior normal to $\Omega $. By Green's Theorem,
we have
\begin{eqnarray*}
\quad
\mu ^+ (\Omega )  & = &
\int\limits_{\partial \Omega } |x|^{1-n} \mbox{ exp } \{ a(|x|)  \} \,
{\mathcal H} ^{n-1} (dx)
\\
     & \geq &
\int\limits_{\partial \Omega } |x|^{-n} \mbox{ exp } \{ a(|x|)  \}
(x\cdot \nu )\,
{\mathcal H} ^{n-1} (dx)
\\
     & = &
n\omega _n e ^{a(0)} +
\int\limits_{ \Omega }
\mbox{ div }
\Big( |x|^{-n} \mbox{ exp } \{ a(|x|)  \}
x \Big)
\, dx
\\
     & = & n\omega _n e ^{a(0)} +
\int\limits_{\Omega } f(x) \, d\mu
\ \ \geq \ \
     n\omega _n e ^{a(0)} +
\int\limits_{B_R } f(x) \, d\mu
\\
     & = &
\int\limits_{\partial B_R  } |x|^{-n} \mbox{ exp } \{ a(|x|)  \}
(x\cdot \nu )\,
{\mathcal H} ^{n-1} (dx)
= \mu ^+ ( B_R ). \qquad \qquad   \qquad \quad \Box
\end{eqnarray*}
Analogously to Theorem \ref{Polya}, the following result holds
\begin{myth}
Let $u$ a nonnegative Lipschitz continuous function with compact
support, and
suppose that
\begin{equation}
\label{6.100}
u(0)= \mbox{ess sup }u.
\end{equation}
Then (\ref{5.1}) holds.
\end{myth}
\noindent{\bf Proof :\/} The proof can be carried out analogously as the proof of
Theorem 6.1 in
\cite{[BBMP2]}, taking into account that the superlevel
sets
$\{ x:\, u(x)>t\} $, ($\sup u >t>0$),
contain an open neighborhood of the origin.
We leave the details to the reader.
$\hfill \Box $

\section{Comparison results}

In this Section we assume that $\mu$ is the measure defined by \eqref{4.800}
with $\varphi$ given by \eqref{4.1999}.

We will prove a comparison result for weak solutions $u$ to nonlinear
elliptic problems. This implies an
estimate of the  Schwarz $\mu $-symmetrization in terms of
the solution of a related radially symmetric problem. We mention that similar
results for the classical
Schwarz symmetrization are well-known (see, for instance,
\cite{[ALT]}, \cite{[Ta1]}, \cite{AT}).
We also mention that a related comparison theorem holds
for the symmetrization in Gauss space (see \cite{[BBMP1]}, \cite{DiBF}).
\bigskip

Consider the following  nonlinear elliptic problem
\begin{equation}
\left\{\begin{array}{ll}
-\hbox{div}(a(x,u,\nabla u))=f\varphi& \hbox{in }
\Omega
\\
u=0 & \hbox{on } \partial\Omega.\label{eq1}
\end{array}
\right.
\end{equation}
Here $\Omega$ is an
open subset of $\R^N$, $N \geq 2$,
$p$ is a real number with $1<p<N$, and
    $a: \Omega \times \R \times \R^N \rightarrow \R^N$
is a Carath\'eodory function
satisfying
\begin{equation}
a(x,s,\xi)\xi \geq \varphi(|x|)|\xi|^p, \label{ip2}
\end{equation}
\begin{equation}
|a(x,s,\xi)| \leq \varphi(|x|) \left[|\xi|^{p-1}+|s|^{p-1} + a_0(x)
\right],
\quad a_0(x)\in L^{p'}(\mu, \Omega),
\ a_0 \geq 0,
\label{ip1}
\end{equation}
\begin{equation}
\left( a(x,s,\xi)-a(x,s,\eta), \xi-\eta \right)>0, \quad \xi
\not=\eta,\label{mon}
\end{equation}
for almost every
$x \in \Omega$ and for every $s \in \R$, $\xi \in \R^N$, $\eta  \in
\R^N$.
Moreover we assume that
\begin{equation}
f\in L^{{(p^*)}'} ( \varphi, \Omega).
\end{equation}
We will say that $u\in   W^{1,p}_0(\Omega, \varphi)$ is a weak solution
to the problem
(\ref{eq1}) if it satisfies
\begin{equation}\label{ws}
\int_\Omega a(x,u,\nabla u)\nabla \psi\, dx=\int_\Omega f\psi\varphi\,
dx, \qquad
\forall \psi \in W^{1,p}_0(\Omega,  \varphi).
\end{equation}
The existence of  a weak solution is a consequence of the Sobolev
type inequality given by Theorem \ref{sobolev} and an adaptation of
classical results due to  J. Leray and J.-L. Lions (cf.
\cite{[LL]}).

The main result of this Section is

\begin{myth}\label{comp2}
Let $u$  be a weak solution to the problem
(\ref{eq1}). Denote by $v=v^{\displaystyle\star}(|x|) \in W^{1,p}_0
(\Omega^{\displaystyle\star}, \varphi)$  the
function defined by
$$
v(|x|)=\int_{H(|x|)}^{\mu(\Omega^{\displaystyle\star})}
\frac{1}{\left[I (r)\right]^{p^\prime}}
\left (\int\limits_{0}^{r}
\tilde f(\sigma) \, d\sigma\right)^{\frac{1}{p-1}}\, dr,
$$
which is a weak solution to the problem
\begin{equation}\label{7.13}
\left\{
\begin{array}{lll}
-\hbox{div}(\varphi (x) |\nabla u|^{p-2}\nabla u)
= f^{\displaystyle\star}(x)  \varphi (x)  & \hbox{in }
     \Omega^{\displaystyle\star},\\
     v=0   &\hbox{on }   \partial\Omega^{\displaystyle \star}.
\end{array}
\right.
\end{equation}
Then
\begin{equation}
\label{7.14}
u^{\displaystyle\star} (x) \leq
    v^{\displaystyle\star} (x), \quad
\end{equation}
for a.e. $x \in   \Omega ^{\displaystyle\star}$. Moreover, for any $1
\le q<p$, there holds
\begin{equation}
\label{grad}
\|\nabla u\|_q \leq
\|\nabla v\|_q.
   \end{equation}
\end{myth}

\bigskip

\noindent{\bf Proof :\/} Denote by $t\in [0, \esup |u|[$, $h>0$ and $\psi_h$
the function
defined by
$$
\psi_h= \left\{
\begin{array}{ll}
\sign u & \qquad \mbox{if } \ |u|> t+h \\
\disp {u-t\sign u\over h } & \qquad \mbox{if } \ t< |u |\le t+h \\
0 & \qquad   \mbox{otherwise. }
\end{array}
\right.
$$
Since $u$ belongs to  $W^{1,p}_0(\Omega, \varphi)$ the function $\psi_h$
can be
considered as a test function in (\ref{ws}) and we have
$$
\int_\Omega a(x,u,\nabla u)\nabla \psi_h\, dx=\int_\Omega f\psi_h
\varphi \, dx,
$$
or, equivalently,
$$
{1\over h}\int\limits_{ t<|u|\le
t+h}a(x,u,\nabla u)\nabla u \, dx =
\int\limits_{|u|>t+h}{f\sign u } \varphi \, dx \ + \ {1\over h}
\int\limits_{t<|u|\le t+h}{f(u-t\sign u)} \varphi
\,dx.
$$
By using the ellipticity condition (\ref{ip2}), Hardy's inequality, and by
letting  $h$ go
to zero, we have that
\begin{eqnarray}\label{5.5}
-\dt \int\limits_{|u|> t}\varphi (|x|)|\D u|^p \, dx
\le \int\limits_{0}^{m _u (t)}
\tilde f(\sigma) \, d\sigma .
\end{eqnarray}
Moreover, by H\"older inequality, we find
\begin{equation}
-\dt \int\limits_{|u|> t}\varphi (|x|)|\D u| \ dx\le  \left (
-\dt \int\limits_{|u|> t} \varphi (|x|)|\D u|^p \ dx\right
)^{1/p} ( -m_u  '(t)  )^{1/p'}
.
\label{5.6}
\end{equation}
On the
other hand, from Federer's coarea formula (cf. \cite{[Fe]}), we obtain
\begin{equation}
-\dt \int\limits_{|u|> t}  {\varphi (|x|)} |\D u| \, dx=
\int\limits_{|u|= t}  {\varphi (|x|)}  \, {\mathcal H}_{n-1} (dx).
\label{5.7}
\end{equation}
Combining (\ref{5.5}), (\ref{5.6}) and (\ref{5.7}), we deduce
\begin{equation}\label{new}
\left(\int\limits_{|u|= t} \varphi (|x|)  \, {\mathcal H}_{n-1} (dx)
\right )^p\le
(-m_u '(t))^{p/p^\prime} \int\limits_{0}^{m _u (t)}
\tilde f(\sigma) \, d\sigma .
\end{equation}
Now we apply the isoperimetric inequality given by Corollary
\ref{isosf}, that is
\begin{equation}\label{disiso}
\int\limits_{u = t} \varphi (|x|)  \, {\mathcal H}_{n-1} (dx) \ge h(H^{-1}
(m_u (t)))=I \bigl(m_u(t)\bigr).
\end{equation}
Combining (\ref{new}) and (\ref{disiso}), we get
\begin{equation}
    -\frac{\left[I \bigl(m_u(t)\bigr) \right]^p}{(-m_u'(t))^{p/p^\prime}}\le
\int\limits_{0}^{m _u (t)}
\tilde f(\sigma) \, d\sigma,  \qquad t\ge \einf u.
\label{A}
\end {equation}
This implies
$$
-\frac{d \tilde u}{ds}\le
\frac{1}{\left[I(s))\right]^{p\prime}}
\Biggl( \int\limits_{0}^{s} \tilde f(r) \,
dr\Biggr)^{\frac{1}{p-1}}.
$$
Now an integration between $s>0$ and $\mu(\Omega^{\displaystyle
\star} )$
gives
$$
\tilde u(s)\le\int_s^{\mu(\Omega^{\displaystyle\star})}
\frac{1}{\left[I(r))\right]^{p^\prime}}
\left (\int\limits_{0}^{r}
\tilde f(\sigma) \, d\sigma\right )^{\frac{1}{p-1}}\, dr.
$$
Choosing $s=H(|x|)>0$ this concludes the proof of (\ref{7.14}).

By H\"older's inequality, we have for any $1\le q<p$,
\begin{equation}
-\dt \int\limits_{|u|> t} \varphi (|x|)  |\D u|^q \, dx\leq
\left (-\dt\int\limits_{|u|> t}   \varphi (|x|)  |\D u|^p\,  dx \right
)^{q/p}(-\mu _u '(t))^{1-q/p}.
    \end{equation}
Using (\ref{5.5}) this leads to
\begin{equation}
-\dt \int\limits_{|u|> t} \varphi (|x|) |\D u|^q \, dx\leq
\left (\int\limits_0^{\mu_u(t)}  \tilde f(s) ds \right
)^{q/p}(-\mu _u '(t))^{1-q/p} .
    \end{equation}
Integrating   between $0$ and $+\infty$ then gives
\begin{equation}
\int\limits_{\Omega} \varphi (|x|)  |\D u|^q \, dx\leq
\int\limits_0^{+\infty} {1\over [ -\mu'(t)]^{q/p}}
\left (\int\limits_0^{\mu_u(t)} \tilde  f(s) ds\right)^{q/p} (-d\mu (t)),
    \end{equation}
from which one has, by (\ref{A}),
\begin{equation}
\int\limits_{\Omega} \varphi (|x|)  |\D u|^q \, dx\leq
   \int\limits_0^{+\infty}
  \left( \frac{1}{I(s)}\right)^{q/(p-1)} \left
(\int\limits_0^s   f^*(r) dr \right )^{q\over p(p-1)}ds.
    \end{equation}
This is (\ref{grad}).
$\hfill \Box $
\medskip

\begin{myrem}\rm We emphasize that the proof of the
comparison result carries over to domains with infinite $\mu$- measure.
Indeed, since
the solutions $u$ and $v$ belong to weighted Sobolev spaces, their
level sets $\{x\in \Omega:˜|u|>t \}$ and $\{x\in \Omega:˜|v|>t \}$
have finite measure and therefore we can apply the isoperimetric
inequality (\ref{4.1102}) to such sets. Notice that in this case
one has to replace $\Omega^{\displaystyle\star}$ in the
symmetrized problem (\ref{7.13}) by $\R^n$.
\end{myrem}

\vspace*{1cm} {\bf authors addresses:}
\\[0.2cm]
Friedemann Brock\\
American University of Beirut\\
 Department of Mathematics\\
 Bliss Street\\
 P.O. Box 11-0236\\
  Beirut\\
   Lebanon\\
   e-mail: fb@13aub.edu.lb
   \\[0.2cm]
Anna Mercaldo \\
Dipartimento
      di Matematica e Applicazioni ``R. Caccioppoli''\\
       Universit\`a
degli Studi di Napoli ``Federico II"\\
 Complesso Monte S. Angelo\\
via Cintia\\
 80126 Napoli\\
  Italy\\
   e-mail:  mercaldo@unina.it
   \\[0.2cm]
   Maria Rosaria Postarero\\
Dipartimento
      di Matematica e Applicazioni ``R. Caccioppoli''\\
       Universit\`a
degli Studi di Napoli ``Federico II"\\
 Complesso Monte S. Angelo\\
via Cintia\\
 80126 Napoli\\
  Italy\\
email: posterar@unina.it

\end{document}